\documentclass[11pt,a4paper, twoside]{article}
\usepackage[notref,notcite]{showkeys}%to show label in pdf 

\usepackage{amsmath,amssymb,amsthm,amsfonts,mathrsfs,amscd,environ}
\usepackage{dsfont}
\usepackage{latexsym,enumerate,color,geometry,extarrows}
\usepackage{xcolor}
\usepackage{verbatim,fancyhdr,enumitem}%for comment paragraph fancyhdr,
 \NewEnviron{ews}{%
\begin{equation}\begin{split}
  \BODY
\end{split}\end{equation}
}

\NewEnviron{ews*}{%
\begin{equation*}\begin{split}
  \BODY
\end{split}\end{equation*}
}

\def\beg{\begin}
\def\bequ{\begin{equation}}
\def\enqu{\end{equation}}
\def\bes{\begin{split}}
\def\ens{\end{split}}
\def\bews{\begin{ews}}
\def\beqn{\begin{eqnarray}}
\def\enqn{\end{eqnarray}}
\def\beq*{\begin{equation*}}
\def\enq*{\end{equation*}}
\def\bqn*{\begin{eqnarray*}}
\def\eqn*{\end{eqnarray*}}
\def\bary{\begin{array}}
\def\eary{\end{array}}
\def\bpma{\begin{pmatrix}}
\def\epma{\end{pmatrix}}
\def\bvma{\begin{Vmatrix}}
\def\evma{\end{Vmatrix}}

 \numberwithin{equation}{section}

\def\al{\alpha}
\def\be{\beta}
\def\ga{\gamma}
\def\de{\delta}
\def\ep{\epsilon}

\def\th{\theta}

\def\la{\lambda}
\def\rh{\rho}

\def\si{\sigma}
\def\ta{\tau}
\def\ph{\phi}

\def\ps{\psi}

\def\Ga{\Gamma}
\def\Us{\Upsilon}
\def\De{\Delta}

\def\Om{\Omega}

\def\R{\mathbb R}
\def\P{\mathbb P}
\def\E{\mathbb E}
\def\D{\mathbb D}

\def\N{\mathbb N}

\def\M{\mathbb M}

\def\sE{\mathscr E}
\def\sF{\mathscr F}

\def\sH{\mathscr H}

\def\cH{\mathcal H}

\def\d{\mathrm{d}}

\def\ff{\frac}
\def\ra{\rightarrow}
\def\nn{\nabla}

\def\<{\langle}
\def\>{\rangle}
\def\sq{\sqrt}
\def\tld{\tilde}
\def\we{\wedge}
\def\1{\mathds{1}}

\def\Span{\mathrm{span}}

\allowdisplaybreaks
\begin{document}

\title{{\bf Stochastic differential equations driven by fractional Brownian motion with locally Lipschitiz drift and their Euler approximation}%\footnote{*.}
}

\author{
{\bf Shao-Qin Zhang$^*$
$\qquad$ 
\bf Chenggui Yuan$^{\dag}$
}\\
~\\
\footnotesize{$\dag$ School of Statistics and Mathematics, Central University of Finance and Economics, Beijing 100081, China.}\\
\footnotesize{ Email: zhangsq@cufe.edu.cn}\\
\footnotesize{$*$ Department of Mathematics, School of Physical Sciences,  Swansea University, Wales SA2 8PP, UK.}\\
\footnotesize{ Email: C.Yuan@swansea.ac.uk}
}

\maketitle

%\linenumbers

\begin{abstract}
In this paper, we study a class of one-dimensional stochastic differential equations driven by fractional Brownian motion with Hurst parameter $H>\ff 1 2$. The drift term of the equation is locally Lipschitz  and unbounded in the neighborhood of $0$. We show the existence,  uniqueness and positivity of the solutions. The estimations of moments, including the negative power moments, are given. Based on these estimations, strong convergence of   the positivity preserving drift-implicit Euler-type scheme is proved, and optimal convergence rate is obtained. By using Lamperti transformation, we show that our results can be  applied to interest rate models such as mean-reverting stochastic volatility model and strongly nonlinear A\"it-Sahalia type model.
\end{abstract}\noindent

\noindent


\vskip 2cm

%\fancyhf{}
%\fancyhead[CO]{Reflection in Hilbert Spaces}
%\fancyhead[CO]{Shao-Qin Zhang}
%\fancyfoot[c]{\thepage}
%\renewcommand{\headrulewidth}{0.4pt}
\section{Introduction}

In this paper, we shall consider a one-dimensional  stochastic differential equation (in short SDEs) driven by fractional Brownian motion: 
\bequ\label{equ-ran}
\d X_t=B(t,X_t)\d t+\si\d B_t^H,~X_0\geq  0,
\enqu
where $B_t^H$ is a fractional Brownian motion (fBM for short) with Hurst $H\in(1/2,1)$ and the drift term $B(t,x)$ is only local Lipschitz in $x\in(0,\infty)$ and unbounded in the neighborhood of $0$. Some general results on this type of equations has been obtained in \cite{HuNS} motivated by the study on Cox-Ingersoll-Ross (C-I-R for short) model in mathematical finance (see \cite{CIR}) with Brownian motion replaced by fBM. Due to the memory effects of fractional Brownian motion, it would be  reasonable to replace Brownian motion by fBM if the there are inert investors in this market, see for instance \cite{RS}. In fact, to handle the complexity of the market, various interest rate models have been developed besides C-I-R, see for instance \cite{Ai-S,BM,CKLS}.  Some of them cannot be covered by the general conditions introduced in \cite{HuNS}. Hence one aim of this paper is to give some general conditions to cover more interest rate models by using the Lamperti transformation, even their coefficients have super-linear growth, see e.g. Example \ref{exam-Ai} the Ait-Sahalia-type interest rate model for details.

On the other hand, numerical approximations of SDEs arising from finance are of great interest. For instance,  strong approximation of C-I-R model based on the Euler-type method was showed in \cite{DNS} and optimal convergence rate is obtained; strong convergence of Euler-Maruyama type approximations for A\"it-Sahalia type model is given \cite{SMHP};  in \cite{MTY}, Euler approximations for  a general  mean-reverting stochastic volatility model under regime-switching is presented. There are many SDEs in mathematical finance having non-Lipschitiz coefficients.  For Euler scheme of SDEs without global  Lipschitz coefficients, one can see \cite{DNS,SMHP,BBD,GR} and references therein. However, the numerical issues for SDEs driven by fBM have not been well studied, comparing with SDEs driven by Brownian motion.  Recently,  the authors in \cite{HHKW}  obtained optimal strong convergence rate of backward Euler scheme for C-I-R model driven by fBM. For numerical scheme of fractional SDEs, one can consult \cite{HHKW,HuLN1, HuLN2,Neu} and references there in.  In this paper, after a general discussion on \eqref{equ-ran}, we  investigate the numerical approximation of the solution to this equation when  $X_0$ is positive. Strong convergence of the numerical scheme is obtained. Based on the Lamperti transformation used as in \cite{DNS,HHKW,HuNS},   our results can cover interesting models in mathematical  finance, such as mean-reverting stochastic volatility model (Example \ref{exam-fra}): 
$$\d Y_t=(a_1-a_2Y_t)\d t+\si Y_t^{\ga}\d B_t^H,~Y_0>0,$$
where $\ga\in [1/2,1)$; and A\"it-Sahalia type model (Example \ref{exam-Ai}): 
$$\d Y_t=\left(a_{-1}Y_t^{-1}-a_0 +a_1Y_t-a_2 Y_t^{r}\right)\d t+\si Y_t^\rh\d B_t^H, Y_0>0,$$
where $\rh>1$ and the stochastic integral in these two models is in the sense of pathwise Riemann-Stieltjes integral developed  by Z\"ahle in \cite{Za}.   The first model was studied in \cite{MTY} under regime-switching, and the convergence rate is obtained. The second model was studied in \cite{SMHP}, where the convergence rate is not clear. 
Following the study in \cite{DNS, HHKW}, the positivity preserving drift-implicit Euler-type method is adopted in our paper. Here, not only  the strong convergence is showed, but also the convergence rate is obtained. For concrete examples presented above, the convergence order of the mean-reverting stochastic volatility model is the Hurst parameter $H$ up to a   logarithmic term, which is an extension of \cite{DNS}; the convergence order of the A\"it-Sahalia type model is $(2H-1)\left(\ff 1 {\rh-1}\we 1\right)$ up to a   logarithmic term. 

This paper is structured as follows. In Section 2, we shall recall some basic facts on fractional Brownian motion. Section 3 is devoted to general discussions on \eqref{equ-ran}, including existence and uniqueness of solutions to the equation; (negative-power) moments and modular of continuity estimations.  In Section 4,  we shall  present our results on the numerical approximations of \eqref{equ-ran} and their applications on concrete examples.

\section{Preliminaries}

We shall recall some basic facts about fractional Brownian motion. For more details, we refer  readers to \cite{BHOZ,Nu,Xiao}.

Let $B^H=\{B_t^H, t\in[0,T]\}$ be a fractional Brownian motion with Hurst parameter $H\in(1/2,1)$ defined on the probability space $(\Om,\sF,\P)$, that is, $B^H$ is a Gaussian process which is centered  with the covariance function
$$
\E\left(B_t^HB_s^H\right)= R_H(t,s)=\frac{1}{2}\left(t^{2H}+s^{2H}-|t-s|^{2H}\right).
$$
For each $t\in[0,T]$, let $\mathcal {F}_t$ be the $\sigma$-algebra generated by the random variables $\{B_s^H:s\in[0,t]\}$ and the sets of
probability zero. Furthermore, one can show that $\E|B_t^H-B_s^H|^p=C(p)|t-s|^{pH}$ for all $p\geq 1$.
As a consequence of the Kolmogorov continuity criterion,
$B^H$ have $(H-\epsilon)$-order H\"{o}lder continuous paths for all $\epsilon>0$. Indeed, the studies on the sample path property of fractional Brownian motion, see for instance \cite{Xiao},  show that 
\beg{align*}
|B_t^H-B_s^H|\leq A|t-s|^H\sq{\log\left(1+(t-s)^{-1}\right)}
\end{align*}
where $A$  is a random variable depending on $H$ only and there is some $c>0$ such that $\E e^{c A^2}<\infty$.

Denote  by $\sE$ the set of step functions on $[0,T]$.
Let $\mathcal {H}$ be the Hilbert space defined as the closure of
$\mathscr{E}$ with respect to the scalar product
$$
\langle I_{[0,t]},I_{[0,s]}\rangle_\mathcal {H}:= \al_H\int_0^T\int_0^T\1_{[0,t]}(u)\1_{[0,s]}(v)|u-v|^{2H-2}\d u\d v=R_H(t,s).
$$
where $\al_H=H(2H-1)$. By the B.L.T. theorem,
the mapping $I_{[0,t]}\mapsto B_t^H$ can be extended to an isometry between $\mathcal {H}$ and the Gaussian space  associated with $B^H$.
Denote this isometry by $\phi\mapsto B^H(\phi)$.

On the other hand, the covariance kernel $R_H(t,s)$ can be written as
$$
 R_H(t,s)=\int_0^{t\wedge s}K_H(t,r)K_H(s,r)\d r,
$$
where $K_H$ is a square integrable kernel given by
\begin{align*}
K_H(t,s)&=\ff {s^{1/2-H}} {\Gamma(H-1/2)}\int_s^t r^{H-1/2}(r-s)^{H-3/2}\d r\1_{[0,t]}(s)
\end{align*}
in which $\Ga(\cdot)$ is the Gamma function. Using this kernel, we could define a map from $L^2([0,T])$ to the reproducing kernel space $\sH$  defined as follows
\beg{align*}
\sH=\overline{\Span\{R_H(t,\cdot)~|~t\in [0,T]\}}^{\<\cdot,\cdot\>_R},\qquad \<R_H(t,\cdot),R_H(s,\cdot)\>_R=R_H(t,s),~s,t\in[0,T].
\end{align*}
For any $\ph\in L^2([0,T])$, let 
$$(K_H\ph)(t)=\int_0^t K_H(t,s)\ph(s)\d s,~t\in [0,T].$$
It has been proved in \cite{BaV,DUs98} that $K_H$ is an isomorphism from $L^2([0,T])$ to $\sH$.

Now, define the linear operator $K_H^*:\mathscr{E}\rightarrow L^2([0,T])$ by
$$
(K_H^*\phi)(s)=K_H(T,s)\phi(s)+\int_s^T(\phi(r)-\phi(s))\frac{\partial K_H}{\partial r}(r,s)\d r.
$$
By integration by parts, it is easy to see that this can be rewritten as
$$
(K_H^*\phi)(s)=\int_s^T\phi(r)\frac{\partial K_H}{\partial r}(r,s)\d r.
$$
It is clear that $\left(K_H^*\1_{[0,t]}\right)(s)=K_H(t,s)\1_{[0,t]}(s)$. $K_H^*$ is the dual operator of $K_H$ in the following sense: for any $\ps\in \sE$ and $h\in L^2([0,T])$, 
\beg{align}\label{dual-K}
\int_0^T (K_H^*\ph)(r)h(t)\d r=\int_0^T \ph(r)(K_H h)(\d r).
\end{align}
Due to \cite{AMN}, for all $\phi,\psi\in\mathscr{E}$,
there holds $\langle K_H^*\phi,K_H^*\psi\rangle_{L^2([0,T])}=\langle\phi,\psi\rangle_\mathcal {H}$
and then $K_H^*$ can be extended to an isometry between $\mathcal{H}$ and $L^2([0,T])$.
Hence, according to \cite{AMN} again,
the process $\{W_t=B^H((K_H^*)^{-1}{\rm I}_{[0,t]}),t\in[0,T]\}$ is a Wiener process,
and $B^H$ has the following integral representation
\beg{align*}
 B^H_t=\int_0^T(K_H^*\1_{[0,t]})(s)\d W_s=\int_0^t K_H(t,s)\d W_s.
\end{align*}

With linear operators $K_H$ and $K_H^*$ in hand, there exists an isometry from $\cH$ to $\sH$ defined by the operator $K_HK_H^*$. Then $\sH$ can be charactered by $\cH$ with  the isometry  $K_HK^*_H$. It follows from the integral representation of fBM that $\sH$ is the fractional version of the Cameron-Martin space. This was showed rigorously in \cite{DUs98}. The Malliavin derivative of the functional of fBM is defined as an $\cH$-valued random variable.   For more details on the Malliavin calculus for fBM, one can consult \cite{Nu}.

In this paper, the stochastic integral of fractional Brownian motion is defined by  the techniques of fractional calculus developed by Z\"{a}hle in \cite{Za}. We cite the following results on the Riemann-Stieltjes integral and chain rule as  a proposition for future use.

\beg{prp}\label{prop-1}
Let $a,b\in\R$ with $a<b$, and let $F\in C^1(\R)$. 
\beg{description}[align=left,noitemsep]
\item {(1)}  Suppose $f\in C^\la(a,b)$ and $g\in C^\mu(a,b)$, where $C^\la(a,b)$ and $C^\mu(a,b)$ are H\"older continuous functions with order $\la$ and $\mu$ respectively. If $\la+\mu>1$, then the Riemann-Stieltjes integral $\int_a^b f\d g$ exists.

\item {(2)} Suppose $f\in C^{\la}(a,b)$ such that $F'\circ f\in C^\mu(a,b)$ with $\la+\mu>1$. Then 
$$F(f(t))-F(f(s))=\int_s^t F'\circ f(r)\d r,~s,t\in (a,b).$$
\end{description}
\end{prp}

Finally, we shall recall a result on the relationship of stochastic integral and the Skorohod integral w.r.t. fractional Brownian motion. Let 
$$|\cH|=\left\{\ps\in \cH~\left|~\|\ps\|_{|\cH|}^2=\al_H\int_0^T\int_0^T|\ps(s)||\ps(t)||t-s|^{2H-2}\d s\d t<\infty \right.\right\}$$ 
and $|\cH|\otimes|\cH|$  be  the set of all measurable function such that
\beg{align*}
\|\ps\|_{|\cH|\otimes|\cH|}^2:=\al_H^2\int_{[0,T]^4}|\ps(u,s)||\ps(v,t)||u-v|^{2H-2}|t-s|^{2H-2}\d u\d v\d t\d s<\infty.
\end{align*}
For $p>1$, we denote by $\D^{1,p}_{|\cH|}$ all the random variable  $u$ such that $u\in |\cH|$  a.s.,   its  Malliavin derivative $Du\in |\cH|\otimes|\cH|$ a.s., and 
$$\E \|u\|_{|\cH|}^p+\E\|Du\|_{|\cH|\otimes|\cH|}^p<\infty.$$
Then we have the following proposition, see \cite[Proposition 5.2.3]{Nu} and \cite{ANu}.
\beg{prp}\label{prop-2}
Let $u_t$ be a stochastic process in $\D^{1,2}_{|\cH|}$ such that a.s.
$$\int_0^T\int_0^T|D_su_t||t-s|^{2H-2}\d t\d s<\infty.$$ 
Then
$$\int_0^T u_t\d B_t^H=\de(u)+\al_H\int_0^T\int_0^T D_su_t|t-s|^{2H-2}\d t\d s.$$
For $p>\ff 1 H$, 
\beg{align*}
\E\left(\sup_{t\in [0,T]}|\de(u\1_{[0,t]})|^p\right)\leq C\left(\E\int_0^T |u_s|^p\d s+\E\left(\int_0^T\int_0^T|D_ru_s|^{\ff 1 H}\d s\d r\right)^{pH}\right).
\end{align*} 
\end{prp}

\section{A study of SDEs driven by fractional Brownian motion}
In this section,  we shall consider \eqref{equ-ran} following \cite{HuNS}.  To get the existence and uniqueness of this equation, we introduce the following assumptions.
\beg{description}[align=left,noitemsep]
\item [(A1)] The drift term $B:[0,\infty)\times (0,\infty)\ra \R$ is continuous and has continuous derivative  w.r.t. the second variable. There exists  $K_t\geq 0$, nondecreasing in $t$,  such that
\beg{align*}
(B(t,x)-B(t,y))(x-y)\leq K_t(x-y)^2,~x,y\in (0,\infty), t\geq 0.
\end{align*}
\item [(A2)] There exist $x_1>0$, $\al>\ff 1 {H}-1$ and $h_1\in C([0,\infty), (0,\infty))$ such that  
\beg{align*}
B(t,x)\geq h_1(t) x^{-\al},~t\geq 0, x\leq x_1.
\end{align*}
\item [(A3)] There is $x_2>0$ and a nonnegative locally bonded function $h_2$ such that 
\beg{align*}
B(t,x)\leq h_2(t)(x+1),~x\geq x_2, t\geq 0.
\end{align*}
\end{description}
%Note that $B(t,\cdot)$ is unbounded in the neighbourhood of $0$. 

The existence and uniqueness of  solutions to \eqref{equ-ran} follows from the existence and uniqueness of the equation below: 
\beg{align}\label{equ-det}
\d X_t=B(t, X_t)\d t+\d w_t,~x_0\geq 0,
\end{align}
where $w\in C^{\be}([0,T],\R)$ for all $T>0$ with $\be\in (\ff 1 2, H)$ such that $\al>\ff 1 {\be}-1$. We say $f$ is a $\be$-H\"older continuous function on $[s,t]$ if  
$$\|f\|_{s,t,\be}:=\sup_{s\leq s'<t'\leq t}\ff {|f(s')-f(t')|} {(t'-s')^\be}<\infty.$$ 
Sometimes, we use $\|\cdot\|_{\be}$ for simplicity's sake. For a continuous function $f$ on $[s,t]$, we define
$$\|f\|_{s,t,\infty}=\sup_{s\leq r\leq t}|f_r|.$$
Our existence and uniqueness theorem for \eqref{equ-det} reads as follows. 
\beg{thm}\label{thm-solution}
Assume that {\bf(A1)}-{\bf (A3)} hold. 
\beg{description}[align=left,noitemsep]
\item {(1)} For all  $X_0>0$,  it holds that the equation \eqref{equ-det} has  a unique solution $X_t$ and 
$$X\in C^\be([0,T],(0,\infty)),~T>0.$$
\item {(2)} For $X_0=0$, if there exists $t^{0}>0$ such that $B(t,\cdot)$ is non-increasing on $(0,x_1)$ for all $0\leq t\leq t^0$, then \eqref{equ-det} has  a unique solution $X_t$ and $X_t\in (0,\infty)$ for all $t>0$.
\end{description} 
\end{thm}
\beg{proof}
We first prove the uniqueness. Let  $X_t^{[1]}$ and $X_t^{[2]}$ be two solutions of equation \eqref{equ-det} with the same initial values, then
$$X_t^{[1]}-X_t^{[2]}=X_s^{[1]}-X_s^{[2]}+\int_s^t \left( B(r,X^{[1]}_r)-B(r,X^{[2]}_r)\right)\d r,~t\geq s.$$
Combining this with {\bf(A1)},  we have
\beg{align*}
\d \left(X_t^{[1]}-X_t^{[2]}\right)^2&=\left( B(t,X^{[1]}_t)-B(t,X^{[2]}_t)\right)\left(X_t^{[1]}-X_t^{[2]}\right)\d t\\
&\leq K_t\left(X_t^{[1]}-X_t^{[2]}\right)^2\d t.
\end{align*}
Thus, it follows from Gronwall's inequality that $X_t^{[1]}-X_t^{[2]}=0$ for all $t\geq 0$.

We assume that $X_0>0$. Since $B:[0,\infty)\times (0,\infty)\ra \R$ is continuous and has continuous derivative  w.r.t. the second variable, it is clear that \eqref{equ-det} has a continuous local solution. Next, we shall prove that $X_t\in (0,\infty)$ for all $t>0$. Let 
$$\ta_0=\inf\{t\geq 0~|~X_t=0 \},\qquad \ta_n=\inf\{t\geq 0~|~X_t\geq n \},~\in\N.$$
We shall prove $\ta_0=\infty$ and $\displaystyle\lim_{n\ra\infty}\ta_n=\infty$.

If $\ta_0<\infty$, then there is $t_0\in (0,\ta_0)$ such that $X_t\leq x_1$ for all $t\in (t_0,\ta_0]$. Since $B(t,x)>0$ for $x\in (0,x_1), t\geq 0$ and
\bequ\label{inequ-ta0}
0=X_{\ta_0}=X_{t}+\int_{t}^{\ta_0}B(s,X_s)\d s+w_{\ta_0}-w_{t_0},
\enqu
it follows that
\bequ\label{inequ-X-w}
X_t\leq |w_{\ta_0}-w_{t}|\leq \|w\|_\be (\ta_0-t)^\be,~t\in (t_0,\ta_0).
\enqu
On the other hand, following from {\bf(A2)}, \eqref{inequ-ta0} and \eqref{inequ-X-w},
\beg{align*}
\|w\|_\be (\ta_0-t)^\be &\geq |w_{\ta_0}-w_{t}|\geq \int_{t}^{\ta_0}B(s,X_s)\d s\\
&\geq \int_t^{\ta_0} h_1(s) X_s^{-\al}\d s\geq \inf_{s\in [t_0,\ta_0]}h_1(s)\int_{t}^{\ta_0}X_s^{-\al}\d s\\
&\geq \ff {\inf_{s\in [t_0,\ta_0]}h_1(s)} {\|w\|_\be^\al} \int_{t}^{\ta_0}\ff 1 {(\ta_0-s)^{\al\be}}\d s \\
&=\ff {(\ta_0-t)^{1-\al\be}\inf_{s\in [t_0,\ta_0]}h_1(s)} {\|w\|^\al_\be},
\end{align*}
this, together with $\al>\ff 1 {\be}-1$, implies that
$$0=\lim_{t\ra\ta_0^-}(\ta_0-t)^{\al\be+\be-1}\geq \lim_{t\ra\ta_0^-}\ff {\inf_{s\in [t_0,\ta_0]}h_1(s)} {\|w\|_\be^{\al+1}}>0.$$
Hence $\ta_0=\infty$. 

If $\ta_\infty:=\displaystyle \lim_{n\ra\infty}\ta_n<\infty$, then either there exists $t_0$ such that $X_{t_0}=x_2+X_0$ and $X_t\geq x_2+X_0$ for all $t\in (t_0,\ta_\infty)$, or for all $n\in \N$ and $\ep>0$ there exits an interval $(t_0,t_1)\subset (\ta_\infty-\ep,\ta_\infty)$ such that $X_{t_0}=x_2+X_0$ and
$$x_2+X_0\leq \inf_{t\in (t_0,t_1)}X_t\leq n\leq  \sup_{t\in (t_0,t_1)}X_t.$$
In both cases, 
\beg{align*}
X_t&=X_{t_0}+\int_{t_0}^t B(s,X_s)\d s+w_{t}-w_{t_0}\\
&\leq x_2+X_0+\int_{t_0}^t h_2(s)(X_s+1)\d s+w_{t}-w_{t_0}\\
&\leq x_2+X_0+\|w\|_{\be}\ta_\infty+\int_0^{\ta_\infty} h_2(s)\d s+\left(\sup_{s\in [0,\ta_\infty]}h_2(s) \right)\int_{t_0}^t X_s\d s.
\end{align*}
where we use {\bf (A3)} in the second inequality. It follows from Gronwall's inequality that for all $t\in (t_0,t_1)$ or $t\in (t_0,\ta_\infty)$
\beg{align*}
X_t&\leq \left(x_2+X_0+\|w\|_{\be}\ta_\infty+\int_0^{\ta_\infty} h_2(s)\d s\right) \exp\left\{(t-t_0)\sup_{s\in [0,\ta_\infty]}h_2(s)\right\}\\
&\leq \left(x_2+X_0+\|w\|_{\be}\ta_\infty+\int_0^{\ta_\infty} h_2(s)\d s\right) \exp\left\{\ta_\infty\sup_{s\in [0,\ta_\infty]}h_2(s)\right\}.
\end{align*}
Taking supremum of the left hand side in the above inequality: for all $t\in (t_0,t_1)$ in the first case or for all $t, \ep, n$ for the second case, the left hand side is infinite but the right hand side is a finite constant. This is a contradiction. Hence, $\ta_\infty=\infty$.

Finally, we shall deal with the case that $X_0=0$. For $n\in\N$, let $X_t^{[n]}$ be the solution of \eqref{equ-det} with $X_0=1/n$. For $n,m\in\N$, $n<m$, let 
$\ta=\inf\{t\geq 0~|~X_t^{[n]}=X_t^{[m]}\}$. By the uniqueness, $X_t^{[n]}=X_t^{[m]}$ for all $t\geq \ta$, or $\ta=\infty$. It is clear that $X_t^{[n]}>X_t^{[m]}$ if $t<\ta$. Thus, the sequence $\{X_t^{[n]}\}_{n\in\N}$ is non-increasing and nonnegative. Let $n_0\in\N$ such that $\ff 1 {n_0}<x_1$, and let $\ta^{(n_0)}=\inf\{t\geq 0~|~X^{[n_0]}_t\geq x_1\}$. Set $X_t=\displaystyle\lim_{n\ra\infty}X_t^{[n]}$. Then
$$X_{t\we\ta^{n_0}}\leq X^{[n]}_{t\we\ta^{(n_0)}}\leq  X^{[n_0]}_{t\we\ta^{(n_0)}}\leq x_1,~n\geq n_0.$$
Because for any $t\in [0,t^0]$, $B(t,x)$ is non-increasing for $x\in (0,x_1)$, the following inequality follows from the monotone convergence theorem
$$\lim_{n\ra\infty}\int_0^{t\we t^0\we \ta^{n_0}} B(s,X^{[n]}_s)\d s=\int_0^{t\we t^0\we \ta^{n_0}} B(s,X_s)\d s.$$
Taking into account that $X^{[n]}_t$ satisfies \eqref{equ-det}, 
$$X_{t\we t^0\we \ta^{n_0}}=\int_0^{t\we t^0\we \ta^{n_0}} B(s,X_s)\d s+w_{t\we t^0\we \ta^{n_0}}-w_0.$$
Moreover, this inequality yields that
$$\int_0^{t\we t^0\we \ta^{n_0}} B(s,X_s)\d s<\infty.$$
Thus, $B(s,X_s)<\infty$ a.e. $s\in [0, t\we t^0\we\ta^{n_0}]$. By {\bf(A2)}, $X_s>0$ a.e. $s\in [0, t\we t^0\we\ta^{n_0}]$.  Starting from any $X_s>0$ with $s\in [0, t\we t^0\we\ta^{n_0})$, there exists unique solution to \eqref{equ-det} which is positive. Thus, $X_s>0$ for all $s\in (0, t\we t^0\we\ta^{n_0}]$. According to the proof above, $\{X_t\}_{t\in [0,t\we t^0\we\ta^{n_0}]}$  can be extended to a solution for all $t>0$ and $X_t>0$ for all $t>0$.

\end{proof}

\beg{rem}
It is clear that $X_t$ is $\be$-H\"older continuous  on $[0,T]$ for all $T>0$ if $X_0>0$. However, we should remark here that, the solution $X_t$ with $X_0=0$ cannot be $\be$-H\"older continuous on interval contains $0$. Otherwise, there is $C>0$ and $t_0>0$ such that $X_t\leq Ct^\be$ for $t\in [0,t_0]$. Letting $\ta_{x_1}=\{t\geq 0~|~X_t\geq x_1\}$, it follows from \eqref{equ-det} and {\bf(A2)} that
\beg{align*}
X_t&=X_0+\int_0^t B(s,X_s)\d s+w_t-w_0\\
&\geq \int_0^t h_1(s) X_s^{-\al}\d s-\|w\|_\be t^\be\\
&\geq \left(\inf_{s\in [0,t]} h_1(s)\right)\int_0^t \ff 1 {C^\al s^{\be\al}}\d s-\|w\|_\be t^\be\\
&= \ff {\left(\inf_{s\in [0,t]} h_1(s)\right) t^{1-\al\be}} {C}-\|w\|_\be t^\be,~t\leq \ta_{x_1}.
\end{align*} 
Then, recalling that $\al>1/\be-1$ implies $1-\al\be-\be<0$,
\beg{align*}
C+\|w\|_\be \geq \varlimsup_{t\ra 0^+} \ff {X_t+\|w\|_\be t^\be} {t^\be}\geq \varlimsup_{t\ra 0^+} \ff {\left(\inf_{s\in [0,t]} h_1(s)\right) t^{1-\al\be-\be}} {C}=\infty.
\end{align*}
\end{rem}

According to this theorem, the  stochastic equation \eqref{equ-ran} has a unique pathwise solution. Next, we shall study the Malliavin differentiablity of $X_t$.

\beg{lem}\label{lem-Mal}
Assume {\bf (A1)}, {\bf (A2)} and {\bf (A3)} hold. Let $X_t$ be  the solution of \eqref{equ-ran}. Then for all $t>0$, $X_t\in \D^{1,2}_{|\cH|}$ with
$$D_s X_t=\si\exp\left\{\int_s^t \nn B\left(r,X_r\right)\d r \right\}\1_{[0,t]}(s),$$
and the law of $X_t$ has density w.r.t. the Lebesgue measure on $\R$.
\end{lem}
The proof  just follows  the line of \cite[Theorem 3.3.]{HuNS}, and the outline of the proof is presented here for the convenience of readers.
\beg{proof}
By {\bf (A1)}, we have
\bequ\label{nnB}
\ff {B(t,x)-B(t,y)} {x-y}\leq K_t,~t>0, x\neq y,
\enqu
which implies the derivative of $B(t,\cdot)$, denoting by $\nn B(t,\cdot)$, is bounded from above by $K_t$. Let $\ep\in (0,1)$, $h\in \cH$ with $h_0=0$  and 
\beg{align*}
X_t^\ep=X_0+\int_0^tB(r,X_r)\d r+\si B^H_t+\si\ep K_HK_H^*h(t).
\end{align*}
Then
\beg{align*}
X^\ep_t-X_t&=\int_0^t \left(B(r,X^\ep_r)-B(r,X_r)\right)\d r+\si\ep K_HK_H^*h(t)\\
&=\int_0^t \nn B\left(r,X_r^\xi\right)(X^\ep_r-X_r)\d r+\si\ep K_HK_H^*h(t),~ t>0,
\end{align*}
where $X_r^\xi=X_r+\xi_s^\ep (X^\ep_r-X_r)$ and $\xi_s^\ep\in (0,1)$ depends on $s$ and $\ep$.  This equality, along with \eqref{dual-K} (see also \cite[Lemma 2.1.9.]{BHOZ}), implies  
\beg{align*}
X^\ep_t-X_t&=\si\ep\int_0^t \exp\left\{\int_s^t \nn B\left(r,X_r^\xi\right)\d r \right\} ( K_HK_H^*h)(\d s)\\
&=\si\ep\int_0^T K_H^*\left(\exp\left\{\int_\cdot^t \nn B\left(r,X_r^\xi\right)\d r \right\}\1_{[0,t]}(\cdot)\right)(s)K_H^* h(s)\d s.
\end{align*}
Since the continuity of $\nn B(t,\cdot)$, \eqref{nnB} and $K_H^* h\in L^2([0,T])$, it follows from the dominated convergence theorem that the limit 
\beg{align*}
\lim_{\ep\ra 0^+}\ff {X^\ep_t-X_t} {\ep}&=\si\int_0^T K_H^*\left(\exp\left\{\int_\cdot^t \nn B\left(r,X_r\right)\d r \right\}\1_{[0,t]}(\cdot)\right)(s)K_H^* h(s)\d s\\
&=\si\left\<\exp\left\{\int_\cdot^t \nn B\left(r,X_r\right)\d r \right\}\1_{[0,t]}(\cdot),h\right\>_\cH
\end{align*}
holds almost sure and in $L^2(\Om)$. Consequently,
$$D X_t=\si\exp\left\{\int_\cdot^t \nn B\left(r,X_r\right)\d r \right\}\1_{[0,t]}(\cdot).$$
It is clear that $\|DX_t\|_{\cH}>0$, and $\E \|D X_t\|_{\cH}^2<\infty$ follows from \eqref{nnB}. Then the existence of density w.r.t. the Lebesgue measure follows from the classical result of Malliavin calculus, see e.g. \cite[Theorem 2.1.2 or Theorem 2.1.3]{Nu}.

\end{proof}

Next, we shall study the moment estimates of solutions to \eqref{equ-ran}. To this end, we introduce the following assumption.
\beg{description}[align=left,noitemsep]
\item [(A2')]  The condition {\bf (A2)} holds.  There exist $\th>0$ and $h_4\in C([0,\infty), (0,\infty))$ such that  
\beg{align}\label{A2'}
B(t,x)\leq h_4(t)(1+x+ x^{-\th}),~t\geq 0, x>0.
\end{align}
\end{description}
It should be noted that $\th\geq \al$ by  {\bf (A2)} and \eqref{A2'}, and {\bf (A2')} implies {\bf (A3)}.  This assumption is used for positive moment estimate. To give the negative moment estimate, we introduce the following
\beg{description}[align=left,noitemsep]
\item [(A3')] there exists a $q>0$ and a locally bounded nonnegative function $h_3$ such that
\bequ\label{A3'}
(B(t,x))^-\leq h_3(t)(1+x^q),~s\geq 0, x>0
\enqu
where $(B(t,x))^-$ denote the negative part of $B(t,x)$.
\end{description}

We first  consider the negative moments for the solution to \eqref{equ-ran}. 

\beg{lem}\label{lem-neg-mom}
Assume {\bf (A1)}-{\bf (A3)} and {\bf(A3')}. Let $X_t$ be  a solution to \eqref{equ-ran} with $X_0>0$. 
\beg{description}[align=left,noitemsep]
\item {(1)}  Suppose $\al=1$. For  $p\geq 1$  with 
\bequ\label{inequ-h1-p-s}
h_1(s)\geq ((p+1)\vee q )Hs^{2H-1}e^{\int_0^sK^+_u\d u},~s\in [0,T],
\enqu
then
$$\sup_{s\in [0,T]}\E  X_s^{-p}<\infty.$$
If \eqref{inequ-h1-p-s} holds with $p$ replaced by $2(p+2)$, then 
$$\E \sup_{s\in [0,T]} X_s^{-p}<\infty.$$
\item {(2)}  Suppose $\al>1$. Then for all $p>0$ and $T\geq 0$,
$$\E \sup_{s\in [0,T]} X_s^{-p}<\infty.$$
\end{description}
\end{lem}
\beg{proof}
We first prove  
\bequ\label{neg-int}
\sup_{s\in [0,T]}\E X_s^{-p}+\int_0^T\E X_s^{-p}\d s<\infty,
\enqu
where for $\al=1$, we impose \eqref{inequ-h1-p-s}. In fact, due to the H\"older inequality, we only need to prove the claim for large $p$. Thus we assume that $p+1\geq q$. Since $X_t$ is $\be$-H\"older continuous for $\be<H$, applying Proposition \ref{prop-1}, Proposition \ref{prop-2} and Lemma \ref{lem-Mal}, we obtain that
\beg{align*}
(X_t+\ep)^{-p}&=(X_0+\ep)^{-p}-p\int_0^t\ff {B(s,X_s)} {(\ep+X_s)^{p+1}}\d s-\si p\int_0^t(\ep+X_s)^{-(p+1)}\d B_s^H\\
&\leq (X_0+\ep)^{-p}-p\int_0^t\ff {B(s,X_s)} {(\ep+X_s)^{p+1}}\d s-\si p\int_0^t(\ep+X_s)^{-(p+1)}\de B^H_s\\
&\qquad +\si p(p+1)\al_H\int_0^t \int_0^s \ff {s^{2H-1}} {(\ep+X_s)^{p+2}}D_r X_s|s-r|^{2H-2}\d r\d s\\
&\leq (X_0+\ep)^{-p}-p\int_0^t\ff {B(s,X_s)X_s-\si^2(p+1)Hs^{2H-1}e^{\int_0^sK_u^+\d u}} {(\ep+X_s)^{p+2}}\d s\\
&\qquad-\si p\int_0^t(\ep+X_s)^{-(p+1)}\de B^H_s.
\end{align*}
Let 
$$\tld x_1=x_1\we 1\we \left(\ff  {e^{-\int_0^TK_u^+\d u}\min_{s\in [0,T]}h_1(s)} {\si^2 (p+1)HT^{2H-1}}\right)^{\ff 1 {\al-1}}\1_{[\al>1]}+(x_1\we 1)\1_{[\al=1]}.$$
Then 
\beg{align*}
-\ff {B(s,x)} {(\ep+x)^{p+2}}&\leq -\1_{[x\leq \tld x_1]}\ff {h_1(s)} {x^\al (\ep+x)^{p+2}}+\1_{[x\geq \tld x_1]} \ff {h_3(s)(1+x^q)} {(\ep+x)^{p+2}}\\
&\leq h_3(s)\left(\ff 1 {\tld x_1^{p+2}}+\ff 1 {\tld x_1^{p+2-q}}\right)
\end{align*}
and
\beg{align*}
&-\ff {B(s,x)x-\si^2(p+1)Hs^{2H-1}e^{\int_0^sK_u^+\d u}} {(\ep+x)^{p+2}}\\
&\qquad\leq -\ff {h_1(s)x^{-(\al-1)}-\si^2(p+1)Hs^{2H-1}e^{\int_0^sK_u^+\d u}} {(\ep+x)^{p+2}}\1_{[x\leq \tld x_1]}\\
&\qquad\qquad +\ff {h_3(s)(1+x^q)x+(p+1)Hs^{2H-1}e^{\int_0^sK_u^+\d u}} {(\ep+x)^{p+2}}\1_{[x\geq \tld x_1 ]}\\
&\qquad\leq -\ff {h_1(s)\tld x_1^{-\al+1}-\si^2(p+1)Hs^{2H-1}e^{\int_0^sK_u^+\d u}} {x^{p+2}}\1_{[x\leq \tld x_1]}\\
&\qquad\qquad +(p+1)Hs^{2H-1}e^{\int_0^sK_u^+\d u} \tld x_1^{-p-2}+h_3(s)\left(\tld x_1^{-(p+1)}+\tld x_1^{-p-1+q}\right).
\end{align*}
Since \eqref{inequ-h1-p-s} and the definition of $\tld x_1$, there exists $C>0$ depending on $\tld x_1,p,q,\si$ such that
\beg{align*}
(X_t+\ep)^{-p}&\leq (X_0+\ep)^{-p}+C\int_0^t (h_3(s)+s^{2H-1})\d s-\si p\int_0^t(\ep+X_s)^{-(p+1)}\de B^H_s.
\end{align*}
Taking expectation and letting $\ep\ra 0$, \eqref{neg-int} is proved. 

If $\al>1$ or \eqref{inequ-h1-p-s} holds with $p$ replaced by $2(P+2)$, then 
$$\sup_{[0,T]} \E X_t^{-2(p+2)} <\infty.$$
Consequently,  
\beg{align*}
&\int_0^T\int_0^T\E (D_s X_t) X_t^{-(p+2)} (D_v X_u) X_u^{-(p+2)}|u-v|^{2H-2}|t-s|^{2H-2}\d u\d v\d s\d t\\
&\qquad \leq \si^2 e^{2\int_0^TK_u^+\d u}\sup_{[0,T]} \E X_t^{-2(p+2)}  \int_0^T\int_0^T|u-v|^{2H-2}|t-s|^{2H-2}\d u\d v\d s\d t\\
&\qquad <\infty.
\end{align*}
Hence $X^{-(p+1)}\in \D^{1,2}_{|\cH|}$. By Proposition \ref{prop-1}, Proposition \ref{prop-2} and Lemma \ref{lem-Mal} again, there is some $C>0$ depending on $\tld x_1,p,q,\si$ such that
\beg{align*}
X_t^{-p}&\leq X_0^{-p}+C\int_0^t (h_3(s)+s^{2H-1})\d s-p\int_0^t X_s^{-(p+1)}\d B_s^H.
\end{align*}
It follows from the maximal inequality of the Skorohod integral  (see e.g. \cite[Page 293]{Nu} or Proposition \ref{prop-2}) that 
\beg{align*}
&\left(\E\sup_{s\in [0,t]}\left|\int_0^t X_s^{-(p+1)}\d B_s^H\right|^2\right)^{\ff 1 2}\\
&\qquad \leq C\left(\int_0^t\E X_s^{-2(p+1)}\d s+\E\int_0^t\left(\int_0^s (p+1)^{\ff 1 H}X_s^{-\ff {p+2} H} |D_rX_s|^{\ff 1 H}\d r\right)^{2H}\d s\right)^{\ff 1 2}\\
&\qquad \leq C_{p,H}(1+t^{H})e^{\int_0^t K_u\d u}\left(\int_0^t\E X_s^{-2(p+1)}\d s\right)^{\ff 1 2}.
\end{align*}
Then
\beg{align*}
\E\sup_{s\in [0,T]}X_s^{-p}&\leq X_0^{-p}+C\int_0^T (h_3(s)+s^{2H-1})\d s\\
&\qquad+C_{p,H}(1+T^{H})e^{\int_0^T K_u^+\d u}\left(\int_0^T\E X_s^{-2(p+1)}\d s\right)^{\ff 1 2},
\end{align*}
which implies the required conclusion.

\end{proof}

If \eqref{A2'} holds with $B(t,x)$ replaced by $|B(t,x)|$, then we can obtain moment estimates of $|X|_{0,T,\infty}$ by applying  \cite[Theorem 3.1]{FanZ} to $X_t^{1+\th}$. However, if \eqref{A3'} holds, that is  we allow that $|B(t,x)|$ has super-linear growth at infinite, then the following lemma can not be covered by \cite{FanZ}. For $g\in C([0,T],\R^d)$, we   denote by $\M_{g,T}(\cdot)$ the modulus of continuity of $g$ on $[0,T]$, i.e.
$$\M_{g,T}(h)=\sup_{0\leq s,t\leq T, |s-t|\leq h}|g_t-g_s|.$$
\beg{lem}\label{lem-psi-mom}
Assume {\bf (A1)}, {\bf (A2')} and {\bf (A3')}. Let $\{X_t\}_{t\geq 0}$ be a solution of \eqref{equ-ran} with $X_0>0$. \\
If $\al>1$, then for any $p>0$ and $T>0$, we have
\bequ\label{inequ-11}
\E \|X\|^p_{0,T,\infty} <\infty,
\enqu
and 
\bequ\label{inequ-22}
\left(\E\M_{X,T}^p+\E\M_{X^{-1},T}^p\right)^{\ff 1 p}\leq C_{p,T}\left(h+h^H\sq{\log(1+1/h)}\right).
\enqu
If $\al=1$, then for $p>0$, there exists $T>0$ such that \eqref{inequ-11} and \eqref{inequ-22} hold.
\end{lem}
\beg{proof}
Suppose $\al>1$. We first prove that 
\bequ\label{EEX}
\E \left(X_t^p+\int_0^t X_s^p\d s\right)<\infty,~t\geq 0, p>0.
\enqu
In fact, by chain rule, Lemma \ref{lem-Mal}, and Proposition \ref{prop-2}, for any $n>0$
\beg{align*}
\ff {nX_t^p} {n+X_t^p}&=\int_0^t\ff {pn^2X_s^{p-1}} {(n+X^p_s)^2}B(s,X_s)\d s+\int_0^t\ff {\si p n^2 X_s^{p-1}} {(n+X_s^p)^2}\d B_s^H\\
&\leq \int_0^t\ff {pn^2h_4(s)X_t^{p-1}(1+X_s+X_s^{-\th})} {(n+X^p_s)^2}\d s+\int_0^t\ff {\si p n^2 X_t^{p-1}} {(n+X_t^p)^2}\de B_t^H\\
&\qquad +\int_0^t\int_0^s\ff {\al \si pn^2 X_s^{p-2}(n(p-1)^+ -(p+1)X_s^p)} {(n+X_s^p)^3}D_rX_s|r-s|^{2H-2}\d r\d s\\
&\leq \int_0^t\left(\ff {2pnh_4(s)X_s^{p}} {n+X^p_s} +ph_4(s)(1+X_s^{-\th})\right)\d s+\int_0^t\ff {\si p n^2 X_t^{p-1}} {(n+X_t^p)^2}\de B_t^H\\
&\qquad +C_He^{\int_0^tK_r\d r}\int_0^t \ff {\si^2 p(p-1)^+n X_s^{p-2}} {n+X_s^p}s^{2H-1}\d s\\
&\leq C_{t,p,K,H,\si}\int_0^t\left(\ff {nX_s^{p}} {n+X^p_s}(h_4(s)+s^{2H-1})\right)\d s+\int_0^t\ff {\si p n^2 X_t^{p-1}} {(n+X_t^p)^2}\de B_t^H \\
&\qquad+C_{t,p,K,H,\si}\int_0^t\left( h_4(s)(1+X_s^{-\th})+(p-1)^+s^{2H-1}X_s^{-2}\right)\d s,
\end{align*} 
where $C_{t,p,K,H,t}$ is locally bounded in $t$. Then it follows from the Gronwall lemma and Lemma \ref{lem-neg-mom} that
\beg{align*}
\E\ff {nX_t^p} {n+X_t^p}\leq e^{C_{t,p,K,H,\si}t}\int_0^t\E\left( h_4(s)(1+X_s^{-\th})+\ff {(p-1)^+ s^{2H-1}} {X_s^{2}}\right)\d s<\infty,
\end{align*}
which implies \eqref{EEX} by letting $n\ra\infty$.

%We shall prove our claim by induction.   Since
%\beg{align*}
%\d X_t=B(t,X_t)\d t+\si \d B_t^H\leq h_4(t)(1+X_t+X_t^{-\th})\d t+\si\d B_t^H,
%\end{align*}
%there is $C>0$ such that 
%$$\sup_{s\in [0,t]}  X_s\leq  X_0+|\si|\sup_{s\in [0,t]}|B_s^H|+C_t\left(1\int_0^t X_s \d s +\int_0^t X_s^{-\th}\d s\right),~t\geq 0.$$
%Then the Gronwall inequality yields that
%\beg{align*}
%\sup_{s\in [0,t]}  X_s&\leq e^{tC_t}\left\{X_0+|\si|\sup_{s\in [0,t]}|B_s^H|+C_t\left(1+\int_0^t X_s^{-\th}\d s\right)\right\}\\
%&\leq e^{tC_t}\left\{X_0+C_t\right\}+|\si|te^{tC_t}\|B^H\|_\infty+C_te^{tC_t}\int_0^t X_s^{-\th}\d s.
%\end{align*}
%It follows from Lemma \ref{lem-neg-mom} that $\E\sup_{s\in [0,t]}  X_s<\infty$ for all $t\geq 0$. Assuming that for $p\geq 2$ there exists 
%\bequ\label{assum-sup}
%\E\sup_{s\in [0,t]}  X^q_s<\infty,~q\leq p-1,t\geq 0,
%\enqu
Next, we shall prove that 
\bequ\label{Esup0}
\E\sup_{s\in [0,t]}  X^{p}_s<\infty,~t\geq 0,p>0.
\enqu
Indeed, by chain rule, \eqref{EEX} and Lemma \ref{lem-Mal}, we have $X_t^{p-1}\in \D^{1,2}_{|\cH|}$ and 
\beg{align}\label{XX-p}
X^p_t&=\int_0^t X_s^{p-1} B(s,X_s)\d s+\si \int_0^t X_s^{p-1}\d B_s^H\nonumber\\
&\leq \int_0^t h_4(s)(X_s^{p-1}+X_s^p+X_s^{p-\th-1})\d s+\si \int_0^t X_s^{p-1}\de B^H_s\nonumber\\
&\qquad +\si (p-1) \al_H\int_0^t \int_0^s X_s^{p-s} D_rX_s |r-s|^{2H-2}\d r\d s\nonumber\\
&\leq C\int_0^t h_4(s)(1+X_s^p)\d s+|\si|\left| \int_0^t X_s^{p-1}\de B^H_s\right|\nonumber\\
&\qquad +C_{H,\si,p}e^{\int_0^t K_r\d r}\int_0^t X_s^{p-2}|s|^{2H-1}\d s.
\end{align}
The maximal inequality of Skorohod integral yields that  the following inequality holds
\beg{align}\label{Esup}
&\left(\E\sup_{s\in [0,t]} \left|\int_0^s X_r^{p-1}\de B^H_r\right|^2\right)^{\ff 1 2}\nonumber\\
&\qquad \leq C\left(\int_0^t\E X_r^{2(p-1)}\d r+\E\int_0^t\left(\int_0^r (p-1)^{\ff 1 H}X_r^{\ff {p-2} H}|D_uX_r|^{\ff 1 H}\d u \right)^{2H}\d r\right)^{\ff 1 2}\nonumber\\
&\qquad \leq  C_{H,p}(1+t^{H})e^{\int_0^tK_r\d r}\left(\int_0^t\E X_r^{2(p-1)}\d r\right)^{\ff 1 2}.
\end{align}
Combining \eqref{XX-p} and \eqref{Esup} with \eqref{EEX}, we get \eqref{Esup0}.

Next, we shall prove the estimates of modulus of continuous. By {\bf (A2')} and {\bf (A3')}, we have
$$|B(s,x)|\leq (h_3\vee h_4)(s)(1+x^q+x^{-\th})\equiv \tld h(s)(1+x^q+x^{-\th}).$$
Then for any $t>s\geq 0$, 
\beg{align}\label{inequ-be-Hold}
|X_t-X_s| \leq &  \int_s^t |B(r,X_r)|\d r+|\si(B_t^H-B_s^H)|\nonumber\\
\leq & \int_s^t \tld h(r)\left(1+X_r^q+X_r^{-\th}\right)\d r+|\si|\M_{B^H,T}((t-s))\nonumber\\
\leq & \left(\sup_{s\leq r\leq t}\tld h(r)\right)\left(1+\|X\|^q_{s,t,\infty}+\|X^{-1}\|_{s,t,\infty}^{\th}\right)(t-s)\nonumber\\
&\quad+|\si|\M_{B^H,T}((t-s)),
\end{align}
which implies for any $p>\ff 1 {1-\be}$,
\beg{align*}
\E\M_{X,T}(h)^p\leq &  C_{T,p}\left(1+\E\|X\|_{0,T,\infty}^{pq}+\E\|X^{-1}\|_{0,T,\infty}^{\th p}\right)h^p+C_p|\si|^p\E(\M_{B^H,T}(h))^p.
\end{align*}
It follows from $\al>1$, Lemma \ref{lem-neg-mom} and the modulus of continuity of $B^H$ (see e.g. \cite[Theorem 4.2]{Xiao} or \cite[Theorem 6.3.3]{MR}) that 
$$\E\M_{X,T}(h)^p\leq C_{p,T}\left\{h^p+h^{pH}\left(\log\left(1+\ff 1 h\right)\right)^{\ff p 2}\right\}.$$

By the H\"older inequality and the following inequality 
\beg{align*}
\sup_{|t-s|\leq h, s,t\leq T}|X^{-1}_t-X^{-1}_s|&\leq \sup_{0\leq s,t\leq T}\left(\ff 1 {X_tX_s}\right)\sup_{|t-s|\leq h, s,t\leq T}|X_t-X_s|\\
&\leq \left(\sup_{0\leq t\leq T}X_t^{-2}\right)\M_{X,T}(h),
\end{align*}
we get the moment estimate of the modulus of continuity of $X^{-1}$.

For $\al=1$, one can repeat the argument for $\al>1$, and takes note that  negative power moments in (1) of Lemma \ref{lem-neg-mom} hold for small $T$ depending on $p$.

\end{proof}

\section{Numerical approximation}

In this section, we shall  consider the numerical approximation of the following equation 
\bequ\label{add-mequ}
\d X_t=B(X_t)\d t+\si\d B^H_t,~X_0>0.
\enqu
The drift term $B(\cdot)$ satisfies  {\bf (A1)}, {\bf (A2')} and {\bf (A3')}, and all these conditions are independent of time. To ensure the positivity of the numerical scheme, we shall use the backward Euler method as in \cite{HHKW}. Moments estimates obtained in the previous section will be used here. 

In addition to  {\bf (A1)}, {\bf (A2')} and {\bf (A3')},  we shall  impose the following assumptions.
\beg{description}[align=left,noitemsep] 
\item [(H1)] There is $h_0>0$ such that the following equation 
$$U(x)+c\equiv B(x)h-x+c=0$$ 
has a unique positive solution for any $c\in \R$ and $0<h<h_0$. 
\item [(H2)] The drift term $B\in C^2(\R)$, and there are nonnegative constants $p_1, p_2$ and $C>0$ such that
\bequ\label{nnB-nn2B}
|\nn B|(x)+|\nn^2 B|(x)\leq C(1+x^{p_1}+x^{-p_2}),~x>0.
\enqu
\end{description}
Sufficient conditions to ensure {\bf (H1)} are that for $0<h< h_0$
$$\lim_{x\ra 0^+}U(x)=\infty,\qquad \lim_{x\ra\infty}\nn U(x)=-\infty$$
and $\nn U(x)<0$ or $\nn^2 U(x)\neq 0$. 

Let $T>0$, $N\in\N$ such that $h:=\ff T N<h_0$, $t_n= nh$, and let  $\De B_{n+1}^H=B^H_{t_{n+1}}-B^H_{t_{n}}$. Since $X_0>0$, we define 
\bequ\label{Xn-she}
X_{n+1}=X_n+B(X_{n+1})h+\si\De B_{n+1}^H,~n\in \N\cup \{0\}.
\enqu
Due to  {\bf (H1)}, the equation \eqref{Xn-she} has a unique positive solution $X_{n+1}$, $n\geq 0$. Let 
$$X_t^h=\ff {t_{n+1}-t} {t_{n+1}-t_n} X_{n}+\ff {t-t_{n}} {t_{n+1}-t_n} X_{n+1},~t_n\leq t\leq t_{n+1}.$$
For a random variable $\xi$, we  denote  $\|\xi\|_p=\left(\E|\xi|^p\right)^{\ff 1 p}$. Our result on numerical approximation of \eqref{add-mequ} reads as follows.

\beg{thm}\label{thm-app}
Assume {\bf (A1)}, {\bf (A2')}, {\bf (A3')}, {\bf (H1)} and {\bf (H2)} hold. Let $h<h_0\we  K^{-1}$, and let $X_{n+1}$ be defined as above.\\
(1) If $\al>1$, then
\beg{align}
\E\sup_{0\leq n\leq N-1 } | X_{t_{n+1}}-X_{n+1} |^p&\leq C_{T,X_0,\th,H,p,B}h^{pH},\label{inequ-n-n-X}\\
\E\sup_{t\in [0,T] }\left| X_{t}-X^h_{t}\right|^p&\leq C_{T,X_0,\th,H,p,B}h^{pH}\left(\log\left(1+\ff 1 h\right) \right)^{p/2}. \label{inequ-h-t-X}
\end{align}
(2)  If $\al=1$, then for $p>0$, there is $T>0$ such that \eqref{inequ-n-n-X} and \eqref{inequ-h-t-X} hold.
\end{thm}
\beg{proof} 
We only prove the claim for $\al>1$. For $\al=1$, the negative power moments estimates hold for $T$ depending on the given $p>0$ (see Lemma \ref{lem-neg-mom}). Then for $T$ small enough, the arguments for $\al>1$ work well in the small interval, and the claim the can be obtained.
   
(1) We first prove \eqref{inequ-n-n-X}. It follows from the definition of $X_{n+1}$ and the mean value theorem that that
\beg{align}\label{Xt-Xn}
X_{t_{n+1}}-X_{n+1}&=X_{t_n}-X_n+\int_{t_n}^{t_{n+1}}B(X_s)\d s-B(X_{n+1})h\nonumber\\
&=X_{t_n}-X_n-\int_{t_n}^{t_{n+1}}\left(B(X_{t_{n+1}})-B(X_s)\right)\d s\nonumber\\
&\qquad+\left(B(X_{t_{n+1}})-B(X_{n+1})\right)h\nonumber\\
&=X_{t_n}-X_n+\nn B(X_{n+1}+\xi_{n+1}(X_{t_{n+1}}-X_{n+1}))h(X_{t_{n+1}}-X_{n+1})\nonumber\\
&\qquad-\int_{t_n}^{t_{n+1}}\left(\int_{s}^{t_{n+1}}\nn B(X_r)B(X_r)\d r+\si\int_{s}^{t_{n+1}}\nn B(X_r)\d B_r^H\right) \d s,
\end{align}
where $\xi_{n+1}\in (0,1)$. By {\bf (A1)},  $\nn B(x)\leq K$ for all $x>0$.  Then, letting 
$$\De_{n+1}=\nn B(X_{n+1}+\xi_{n+1}(X_{t_{n+1}}-X_{n+1})), $$
we have
$$1-\De_{n+1}h\geq 1-Kh>0$$
holds for small $h$. On the other hand, it follows from the Fubini theorem that
\beg{align*}
&\int_{t_n}^{t_{n+1}}\left(\int_{s}^{t_{n+1}}\nn B(X_r)B(X_r)\d r+\si\int_{s}^{t_{n+1}}\nn B(X_r)\d B_r^H\right) \d s\\
&\qquad =\int_{t_n}^{t_{n+1}}\left(r-t_n\right)\nn B(X_r)B(X_r)\d r+\si\int_{t_n}^{t_{n+1}}(r-t_n)\nn B(X_r)\d B_r^H.
\end{align*}
Substituting this into \eqref{Xt-Xn}, letting $\Us_{n+1}=X_{t_{n+1}}-X_{n+1}$ and 
$$Q_{n+1}=-\int_{t_n}^{t_{n+1}}\left(r-t_n\right)\nn B(X_r)B(X_r)\d r-\si\int_{t_n}^{t_{n+1}}(r-t_n)\nn B(X_r)\d B_r^H,$$
we get that
\beg{align*}
\Us_{n+1}=(1-\De_{n+1}h)^{-1}\Us_n+(1-\De_{n+1}h)^{-1}Q_{n+1}.
\end{align*}
Consequently, 
\bequ\label{equ-Us}
\Us_{n+1}=\sum_{i=1}^{n+1}Q_{i}\prod_{k=i}^{n+1}(1-\De_{k}h)^{-1}=:\sum_{i=1}^{n+1}Q_{i}\rh_i.
\enqu

Next, we shall estimate the right hand side of the above equality. Since
\bequ\label{Inequ-de-1}
\prod_{k=i}^{n+1}(1-\De_{k}h)^{-1}\leq (1-Kh)^{-n+i}\leq e^{n\log \ff 1 {1-Kh}}\leq e^{\ff {nKh} {1-Kh}}=e^{\ff {KT} {1-Kh}},
\enqu
it follows from \eqref{equ-Us} that 
$$\E\sup_{1\leq n\leq N}|\Us_{n}|^p\leq C_{T,K}\E\left( \sum_{i=1}^N |Q_i|\right)^{p}.$$
By the definition of $Q_i$, there are two integrals to be estimated.  For the ordinary integral, it follows from {\bf (A2')}, {\bf (A3')} and {\bf (H2)} that 
\beg{align}\label{Inequ-nnBB}
\left\|\sum_{i=1}^{N} \int_{t_{i-1}}^{t_i}(t_{i-1}-r)\nn B(X_r)B(X_r)\d r \right\|_p&\leq C_{K,T} h\sum_{i=1}^{N}\int_{t_{i-1}}^{t_i}\left\|\nn B(X_r)B(X_r)\right\|_p\d r\nonumber\\
&\leq C_{T} h\int_{0}^{T}\|1+X^{-(\th+p_1)}_r+X_r^{q\vee (1 +p_2)}\|_p\d r\nonumber\\
&\leq C_{T,\th,X_0,p_1,p_2,q,K} h.
\end{align}
For the stochastic integration, by \cite[Theorem 5.2.3]{Nu}
\beg{align}\label{Inequ-I1-I2}
&\sum_{i=1}^{N}\left|\int_{t_{i-1}}^{t_i}(r-t_{i-1})\nn B(X_r)\d B_r^H\right|\nonumber\\
&\qquad\leq \sum_{i=1}^{N}\left|\int_0^T (t-t_{i-1})\nn B(X_t)\1_{(t_{i-1},t_i]}\de B_t^H\right|\nonumber\\
&\qquad\quad +\sum_{i=1}^{N}\left|\int_{0}^{T}\int_{0}^{T}(r-t_{i-1})\nn^2B(X_r)D_sX_r|s-r|^{2H-2}\1_{(t_{i-1},t_i]}(r)\1_{(0,r]}(s)\d s\d r\right|\nonumber\\
&\qquad=\sum_{i=1}^{N}\left|\int_0^T(r-t_{i-1})\nn B(X_r)\1_{(t_{i-1},t_i]}(r)\de B_r\right|\nonumber\\
&\qquad\quad +\sum_{i=1}^{N}\int_{t_{i-1}}^{t_i}\int_{0}^{r}(r-t_{i-1})\left|\nn^2B(X_r)\right||D_sX_r||s-r|^{2H-2}\d s\d r\nonumber\\
&\qquad =: I_1+I_2.
\end{align}
For $I_2$, it follows from \eqref{Inequ-de-1} that
\beg{align}\label{Inequ-I2}
\|I_2\|_p&\leq \left\|\sum_{i=1}^{N}\int_{t_{i-1}}^{t_i}\int_{0}^{r}(r-t_{i-1})|\nn^2B(X_r)||D_sX_r||s-r|^{2H-2}\d s\d r\right\|_p\nonumber\\
&\leq C h\sum_{i=1}^{N}\int_{t_{i-1}}^{t_i}\int_{0}^{r} \|X_r^{-p_1 }+X_r^{p_2}\|_p e^{\int_s^r K\d u}|s-r|^{2H-2}\d s\d r\nonumber\\
%& \leq C_{T,\th,a,b^-} h\sum_{i=0}^{n+1}\int_{t_{i-1}}^{t_i}\int_{t_{i-1}}^{t_i} h_2(r)X_r^{-(\th+2)}e^{\int_s^r b^-\d u}|s-r|^{2H-2}\d s\d r\\
& \leq C_{T,K,H} h\sum_{i=1}^{N}\int_{t_{i-1}}^{t_i}  \|X_r^{-p_1 }+X_r^{p_2}\|_p r^{2H-1}\d s\d r\nonumber\\
& \leq C_{T,\th,H,X_0,K} h.
\end{align}
For $I_1$, it follows from Minkowski's inequality that 
\beg{align*}
\|I_1\|_p&=\left\|\sum_{i=1}^{N}\left|\int_0^T(r-t_{i-1})\nn B(X_r)\1_{(t_{i-1},t_i]}\de B_r^H\right|\right\|_p\\
&\leq\sum_{i=1}^{N}\left\|\int_0^T(r-t_{i-1})\nn B(X_r)\1_{(t_{i-1},t_i]}\de B_r^H\right\|_p.
\end{align*}
By  \cite[Proposition 1.5.8]{Nu}, 
\beg{align*}
&\E \left|\int_0^T(r-t_{i-1})\nn B(X_r)\1_{(t_{i-1},t_i]}\de B_r^H\right|^p\\
&\qquad\leq C_p\left(\int_{(t_{i-1},t_i]^2}(r-t_{i-1})(s-t_{i-1})|\E\nn B(X_r)||\E \nn B(X_s)||r-s|^{2H-2}\d r\d s\right)^{\ff p 2}\\
&\qquad\qquad  +\E\Big(\int_{t_{i-1}}^{t_i}\int_{t_{i-1}}^{t_i}\int_0^r\int_0^s(r-t_{i-1})(s-t_{i-1})|\nn^2B(X_r)||\nn^2 B(X_s)|\\
&\qquad\qquad\qquad \times|D_u X_r||D_v X_s||u-v|^{2H-2}|s-r|^{2H-2}\d u\d v\d s\d r\Big)^{\ff p 2}\\
&\qquad \leq  C_{T,p_1,p_2,K}h^p \left(\int_{(t_{i-1},t_i]^2}|r-s|^{2H-2}\d r\d s\right)^{\ff p 2}\\
&\qquad\qquad +C_{T,K,H} h^p \Big(\int_{(t_{i-1},t_i]^2}\int_0^r\int_0^s\left(\E\left(1+X_r^{-p_1}+X_r^{p_2}\right)^{\ff p 2}\left(1+X_s^{-p_1}+X_s^{p_2}\right)^{\ff p 2}\right)^{\ff 2 p}\\
&\qquad\qquad\qquad \times |u-v|^{2H-2}|r-s|^{2H-2}\d u\d v\d r\d s\Big)^{\ff  p 2 }\\
&\qquad \leq  C_{T,X_0,H,K}h^{p+Hp}\\
&\qquad\qquad+C_{T,X_0,K,H,p_1,p_2,p}h^{p}\left(\int_{(t_{i-1},{t_i}]^2}\int_{(0,T]^2}|u-v|^{2H-2}|r-s|^{2H-2}\d u\d v\d r\d s\right)^{\ff p 2}\\
&\qquad \leq  C_{T,X_0,K,H,p_1,p_2,p}h^{p+pH}.
\end{align*}
Thus 
\beg{align}\label{Inequ-I1}
\|I_1\|_p & \leq C_{T,X_0,p,H,K}\sum_{i=1}^{N} h^{1+H}\leq C_{T,X_0,p,H,K} h^H.
\end{align}
Substituting \eqref{Inequ-nnBB}, \eqref{Inequ-I1-I2}, \eqref{Inequ-I2} and \eqref{Inequ-I1} into \eqref{equ-Us}, we obtain
$$\E\sup_{0\leq n\leq N-1}|\Us_{n+1}|^p\leq C_{T,X_0,\th,H,p,q,K,p_1,p_2} h^{pH}.$$

(2)~
For $t\in [t_n,t_{n+1}]$,
\beg{align*}
\left|X_t-X^h_t\right| &=  \Big|-\ff {t-t_n} {h}X_{t_{n+1}}-\ff {t_{n+1}-t} {h}X_{t_{n}}+X_t+\ff {t-t_n} {h}(X_{t_{n+1}}-X_{n+1})\\
&\qquad +\ff {t_{n+1}-t} {h}(X_{t_{n}}-X_n)\Big|\\
&\leq  \left|X_{t_{n+1}}-X_{t}\right|+\left|X_{t_{n}}-X_{t}\right|+|\Us_{n+1}|+|\Us_n|\\
&\leq  2\int_{t_n}^{t_{n+1}}|B(X_r)|\d r+|B^H_{t_{n+1}}-B^H_t|+|B^H_t-B^H_{t_{n}}|\\
&\qquad +|\Us_{n+1}|+|\Us_n|\\
&\leq  C\int_{t_n}^{t_{n+1}}\left(X_r^{q\vee 1}+X_r^{-\th}\right)\d r+2\M_{H,T}(h)+|\Us_{n+1}|+|\Us_n|\\
&\leq C\left(\|X\|^{q\vee1}_{0,T,\infty}h+\left(\int_{0}^{T}X_r^{-\ff {\th} {1-H}}\d r\right)^{1-H}h^H\right)\\
&\qquad +2\M_{B^H,T}(h)+2\sup_{1\leq n\leq N}|\Us_n|\\
&=: I_1+I_2+I_3.
\end{align*}
Since $I_1$, $I_2$ and $I_3$ are independent of $n$ and $t$, we have
$$\E\sup_{t\in [0,T]}\left|X_t-X^h_t\right|^p\leq 3^{p-1}\left(\E I_1^p+\E I_2^p+\E I_3^p\right).$$
It follows from Lemma \ref{lem-neg-mom} and Lemma \ref{lem-psi-mom} that 
$$\E I_1^p\leq C_{K,q,T,X_0,p,\th,H} h^{pH}.$$
The inequality \eqref{inequ-n-n-X} yields that 
$$\E I_3^p=2^p\E\sup_{1\leq n\leq N}|\Us_n|^p\leq C_{K,q,T,X_0,p,\th,H}h^{pH}.$$
The modulus of continuity of $B^H$ (see e.g. \cite[Theorem 4.2]{Xiao} or \cite[Theorem 6.3.3]{MR}) implies that there is a constant $C_{T,p}>0$ such that
\beg{align*}
\E I_2^p=2^p\E \M_{B^H,T}^p(h)&\leq  C_{T,p}  h^{pH}\left(\log(1+ 1/h) \right)^{p/2},
\end{align*}
Therefore, 
$$\E \sup_{t\in [0,T]}\left|X_t-X^h_t\right|^p\leq C_{K,q,T,X_0,p,\th,H}h^{pH}\left(\log(1+ 1/h) \right)^{p/2}.$$

\end{proof}

Some  concrete models can be transformed  to \eqref{add-mequ}, see  Example \ref{exam-fra} and Example \ref{exam-Ai} for instance. The following corollary is crucial to getting the numerical approximation of them.   
\beg{cor}\label{cor-app-l}
Assume  the hypotheses  of Theorem \ref{thm-app} hold.\\
(1)~If $\al>1$, then for any $l>0$, 
\bequ\label{l-l-1}
\left(\E\sup_{t\in [0,T]}|X_t^l-(X_t^h)^l|^p\right)^{\ff 1 p}\leq Ch^{H(l\we 1)}\left(\log(1+ 1/h)\right)^{\ff {l\we 1} 2};
\enqu
for any $l\in(0,\al]$,
\bequ\label{l-l-2}
\left(\E\sup_{t\in [0,T]}|X_t^{-l}-(X_{t}^h)^{-l}|^p\right)^{\ff 1 p}\leq C_{p,T}h^{(2H-1)(l\we 1)}\left(\log(1+1/h)\right)^{l\we 1};
\enqu
(2)~If $\al=1$, then for $l>0$ and $p>0$, there is $T>0$ such that \eqref{l-l-1} holds; for $l\in (0,1]$ and $p>0$, there is $T>0$ such that \eqref{l-l-2} holds. 
\end{cor}

\beg{proof}
For $l\in(0,1]$, it follows from Lemma \ref{thm-app} that
\beg{align*}
\left(\E\sup_{t\in [0,T]}|X_t^l-(X_t^h)^l|^p\right)^{\ff 1 p}&\leq \left(\E\sup_{t\in [0,T]}|X_t-X_t^h|^{lp}\right)^{\ff 1 p}\leq \left(\E\sup_{t\in [0,T]}|X_t-X_t^h|^p\right)^{\ff l p}\\
&\leq Ch^{lH}\left(\log(1+\ff 1 h)\right)^{\ff l 2}.
\end{align*}
For $l>1$, 
\beg{align*}
&\left(\E\sup_{t\in [0,T]}|X_t^l-(X_t^h)^l|^p\right)^{\ff 1 p}\\
&\qquad\leq \left(\E\sup_{t\in [0,T]}\left(X_t^{p(l-1)}\vee(X_t^h)^{p(l-1)}\right)|X_t-X_t^h|^p\right)^{\ff 1 p}\\
&\qquad\leq\left(\E\sup_{t\in [0,T]}\left(X_t^{2p(l-1)}\vee(X_t^h)^{2p(l-1)}\right)\right)^{\ff 1 2}\left(\E\sup_{t\in [0,T]}|X_t-X_t^h|^{2p}\right)^{\ff 1 {2p}}\\
&\qquad\leq C_{p,T}h^H\left(\log(1+\ff 1 h)\right)^{\ff 1 2}.
\end{align*}
Hence, we have proved our first claim.

To consider the negative power approximation, we first give an estimate of $X_n^{-1}$. By \eqref{Xn-she}, there is positive constant $C$ which is independent of $n,h$ such that
\beg{align*}
&C\left(X_{n+1}^{-\al}h-(X_{n+1}^q+1)h\right)\leq B(X_{n+1})h\\
&\qquad\leq |X_{n+1}-X_{t_{n+1}}|+|X_{t_{n+1}}-X_{t_n}|+|X_{n}-X_{t_n}|+|\si||B^H_{n+1}|.
\end{align*}
Then 
\beg{align}\label{add-ineq}
C\sup_{1\leq n\leq N} X_{n}^{-\al}&\leq C\sup_{1\leq n\leq N}(X_{n}^q+1)+\ff 2 h\sup_{1\leq n\leq N}|X_{n}-X_{t_n}| \nonumber\\
&\qquad+\ff 1 h {\sup_{1\leq n\leq N}|X_{t_{n+1}}-X_{t_n}|} +\ff {|\si|} h \M_{B^H,T}^p(h)\nonumber\\
&\leq C\sup_{1\leq n\leq N}(X_{n}^q+1)+\ff 2 h\sup_{1\leq n\leq N}|X_{n}-X_{t_n}| \nonumber\\
&\qquad+\ff 1 h\M_{X,T}(h) +\ff {|\si|} h \M_{B^H,T}^p(h).
\end{align}
Since \eqref{inequ-n-n-X}, it is clear that for all $p>0$, we have $\E \sup_{1\leq n\leq N}\left(X_{n}^{qp}\right)<\infty$. By Theorem \ref{thm-app},  
$$\E\sup_{1\leq n\leq N}|X_{n}-X_{t_n}|^p\leq h^{Hp}\left(\log(1+\ff 1 h)\right)^{\ff p 2}.$$
It follows from Lemma \ref{lem-psi-mom} that 
$$\E\M_{X,T}^p(h)\leq h^{Hp}\left(\log(1+\ff 1 h)\right)^{\ff p 2}.$$
Combining these with \eqref{add-ineq}, we get that
\bequ\label{EXn-p}
\E\sup_{1\leq n\leq N} X_{n}^{-p\al}\leq C\left(1+h^{(H-1)p}\left(\log(1+\ff 1 h)\right)^{\ff p 2}\right).
\enqu
Then for $l\in [1,\al]$, it follows from \eqref{EXn-p} that
\beg{align*}
&\left(\E\sup_{t\in [0,T]}|X_{t}^{-l}-(X_{t}^h)^{-l}|^p\right)^{\ff 1 p}=\left(\E\sup_{t\in [0,T]}\ff {|X_{t}^{l}-(X_{t}^h)^{l}|^p} {X_{t}^{pl} (X_{t}^h)^{pl}}\right)^{\ff 1 p}\\
&\qquad\leq \left(\E\sup_{t\in [0,T]}X_t^{-3pl}\right)^{\ff 1 {3p}}\left(\E\sup_{t\in [0,T]}(X_{t}^{h})^{-3pl}\right)^{\ff 1 {3p}}\left(\E\sup_{t\in [0,T]}|X_{t}^{l}-(X_{t}^h)^{l}|^{3p}\right)^{\ff 1 {3p}}\\
&\qquad= \left(\E\sup_{t\in [0,T]}X_t^{-3pl}\right)^{\ff 1 {3p}}\left(\E\sup_{1\leq n\leq N}X_{n}^{-3pl}\right)^{\ff 1 {3p}}\left(\E\sup_{t\in [0,T]}|X_{t}^{l}-(X_{t}^h)^{l}|^{3p}\right)^{\ff 1 {3p}}\\
&\qquad \leq C_{p,T}\left(1+h^{H-1}\left(\log(1+\ff 1 h)\right)^{\ff 1 2}\right)h^{H}\left(\log(1+\ff 1 h)\right)^{\ff 1 2}\\
&\qquad \leq C_{p,T}h^{2H-1}\log(1+\ff 1 h).
\end{align*}
For $l<1$, 
\beg{align*}
\left(\E\sup_{t\in [0,T]}|X_{t}^{-l}-(X_{t}^h)^{-l}|^p\right)^{\ff 1 p}&\leq \left(\E\sup_{t\in [0,T]}|X_{t}^{-1}-(X_{t}^h)^{-1}|^p\right)^{\ff l p}\\
&\leq C_{p,T,l}h^{(2H-1)l}\left(\log(1+\ff 1 h)\right)^l.
\end{align*}
Combining these two cases together, we prove our second conclusion.

\end{proof}

\beg{rem}
If $\ph$ is a continuous function on $(0,\infty)$ such that 
$$|\ph(x)-\ph(y)|\leq C|x^l-y^l|\qquad\mbox{or}\qquad |\ph(x)-\ph(y)|\leq C|x^{-l}-y^{-l}|$$
for $l$ as in Corollary \ref{cor-app-l} and some $C>0$. Then we can  approximate  $\ph(X_t)$ by $\ph(X_t^h)$. 

For $\al=1$, the convergence of the backward Euler  scheme for C-I-R model driven by fractional Brownian motion has been obtained in \cite{HHKW}. Theorem \ref{thm-app} and Corollary \ref{cor-app-l} can also be applied to C-I-R model, and stronger convergence  can be obtained for some $p>0$ and some small $T>0$ depending on $p$.  To get more specific and sharp dependencies between $p$ and $T$, one can follows  the proof of   \cite[Theorem 4.1]{HHKW} and Theorem \ref{thm-app}.
\end{rem}

Finally, we apply our  results to the two examples  introduced in the introduction.
\beg{exa}\label{exam-fra}
We consider the numerical simulation of the following equation
\beg{align}\label{equ-Y-1}
\d Y_t=(a_1 -a_2 Y_t)\d t+\si Y_t^{\ga}\d B_t^H, X_0>0
\end{align}
with $\ga\in (\ff 1 2,1)$, $a_1>0$, $a_2\in \R$ and $\si\neq 0$. To study this equation, we consider 
$$\d X_t=(1-\ga)\left(a_1X_t^{-\ff {\ga} {1-\ga}}-a_2X_t\right)\d t+\si (1-\ga)\d B_t^H, X_0=Y_0^{1-\ga}.$$
Setting $B(x)=(1-\ga)a_1 x^{-\ff \ga {1-\ga}}-a_2(1-\ga)x$, it is clear that {\bf (A1)}, {\bf (A2')} and {\bf (A3)}  hold  with $K=-a_2$, $\th=\al=-\ff {\ga} {1-\ga}$, and $q=\1_{[a_2>0]}$. Then this equation has a unique solution by applying Theorem \ref{thm-solution}. Moreover, it follows from the chain rule that $Y_t=X_t^{\ff 1 {1-\ga}}$ and \eqref{equ-Y-1} has a uniqueness solution. It is clear that there is $h_0>0$ such that $1+a_2(1-\ga)h>0$ for $h\in(0, h_0)$. Then for all $c\in\R$, the equation
$$B(x)h-x+c\equiv (1-\ga)a_1 x^{-\ff {\ga} {1-\ga}}h-(1+a_2(1-\ga)h)x+c=0$$
has a unique positive solution. It follows from Corollary \ref{cor-app-l} that
$$\left(\E\sup_{t\in [0,T]}|Y_t-(X_t^h)^{\ff 1 {1-\ga}}|^p\right)^{\ff 1 p}\leq Ch^{H }\left(\log(1+ 1/h)\right)^{\ff {1} 2}.$$
\end{exa}

\beg{exa}\label{exam-Ai}
In this example, we investigate the nonlinear A\"it-Sahalia-type interest rate model:
\beg{align}\label{Ait-Sa}
\d Y_t=\left(a_{-1}Y_t^{-1}-a_0 +a_1Y_t-a_2 Y_t^{r}\right)\d t+\si Y_t^\rh\d B_t^H, Y_0>0,
\end{align}
with $r+1>2\rh$ and $r\geq 2\we\rh +1>2$ and $a_i>0$, $i=-1,0,1,2$.  To study \eqref{Ait-Sa}, we consider
\beg{align}\label{Ait-Sa-X}
\d X_t&=(\rh-1)\left(a_2X_t^{-\ff {r-\rh} {\rh-1}}-a_1 X_t+a_0 X_t^{\ff {\rh} {\rh-1}}-a_{-1}X_t^{\ff {\rh+1} {\rh-1}}\right)\d t\nonumber\\
&\qquad +(1-\rh)\si\d B_t^H, X_0=Y_0^{1-\rh}.
\end{align}
Set 
\beg{align*}
B(x)&=(\rh-1)\left(a_2x^{-\ff {r-\rh} {\rh-1}}-a_1 x+a_0 x^{\ff {\rh} {\rh-1}}-a_{-1}x^{\ff {\rh+1} {\rh-1}}\right)\\
&\equiv b_1x^{-\ff {r-\rh} {\rh-1}}-b_2 x+b_3 x^{\ff {\rh} {\rh-1}}-b_4 x^{\ff {\rh+1} {\rh-1}}.
\end{align*} 
Since $\ff {r-\rh} {\rh-1}>1$, it is clear that {\bf (A1)}, {\bf (A2')} and {\bf (A3')}  hold with $\th=\al=\ff {r-\rh} {\rh-1}$, $q=\ff {\rh+1} {\rh-1}$ and some constant $K$. Then this equation  has a unique solution, and so does \eqref{Ait-Sa}. Moreover $Y_t=X_t^{-\ff  1 {\rh-1}}$. It is clear by $\ff {r-\rh} {\rh-1}>1$ and $\ff {\rh+1} {\rh-1}>\ff {\rh} {\rh-1}$ that for $h>0$
$$\lim_{x\ra 0^+}(B(x)h-x)=+\infty,\qquad \lim_{x\ra +\infty}(B(x)h-x)=-\infty.$$
On the other hand,
\beg{align*}
\nn B(x)h-1&=-\ff {b_1(r-\rh)h} {\rh-1} x^{-\ff {r+1} {\rh-1}}-(b_2h+1)\\
&\qquad+\ff {b_3 \rh h} {\rh-1}x^{\ff 1 {\rh-1}}-\ff {b_4(\rh+1)h} {\rh-1}x^{\ff 2 {\rh-1}}.
\end{align*}
Then for  $0<h<\ff {4(\rh-1)b_4(\rh+1)} {b_3^2\rh^2}$, we have
$$\ff {(b_3\rh)^2h^2} {(\rh-1)^2}-4\left(\ff {b_1(r-\rh)h} {\rh-1}x^{\ff {r+1} {\rh-1}}+b_2h+1\right)\ff {b_4(\rh+1)h} {\rh-1}<0,$$
which implies that $\nn B(x)h-1<0$. Consequently, {\bf (H1)} holds.  Hence, Theorem \ref{thm-app} can be applied to \eqref{Ait-Sa-X}.

Since $r+1>2\rh$ and $r>2\we\rh+1$, we have $\ff 1 {\rh-1}\leq \ff {r-\rh} {\rh-1}$. Letting $l=\ff 1 {\rh-1}$ in Corollary \ref{cor-app-l}, we have 
$$\left(\E\sup_{t\in [0,T]}|Y_t-(X_t^h)^{\ff 1 {\rh-1}}|^p\right)^{\ff 1 p}\leq C_{p,T}h^{(2H-1)(\ff 1 {\rh-1}\we 1)}\left(\log(1+1/h)\right)^{\ff 1 {\rh-1}\we 1},$$
which implies that 
$$\lim_{h\ra 0^+}\E\sup_{t\in [0,T]}|Y_t-(X_t^h)^{\ff 1 {\rh-1}}|^p=0.$$
\end{exa}
%\noindent\textbf{Acknowledgements}

%\medskip

\end{document}